# FURTHER GENERALIZATION OF GOLDEN MEAN IN RELATION TO EULER'S "DIVINE" EQUATION

Miloje M. Rakočević

*Department of Chemistry, Faculty of Science, University of Niš,*
*Ćirila i Metodija 2, Serbia(E-mail: m.m.r@eunet.yu)*

**Abstract**.

In the paper a new generalization of the Golden mean, as a further generalization in relation to Stakhov (1989) and to Spinadel (1999), is presented. Also it is first observed that the Euler's "divine equation" $(a + b^n)/n = x$ represents a possible generalization of Golden mean.

## 1. INTRODUCTION

The Golden mean canon (GM) or ratios close to it are found in the linear proportions of masterpieces of architecture, human, animal, and plant bodies. In last decades the canon is extended to the periodic system of chemical elements (Luchinskiy & Trifonov, 1981; Rakočević, 1998a, Djukić & Rakočević, 2002) and to genetic code (Rakočević, 1998b), as well as to the different natural and artificial structures, especially to nanotehnology (Koruga et al., 1993; Matija, 2004). As a noteworty fact, the Golden mean is found in masterpieces of classic literature (Stakhov, 1989; Freitas, 1989; Rakočević, 2000; Rakočević, 2003).

In the present day there are minimum two generalisations of Golden mean. First, a "vertical" generalization with $x^n$ instead $x^2$ in the equation of Golden mean [Equations (2) and (4) in the next Section] (Stakhov, 1989); second, a "horizontal" generalization with $p > 1$ and/or $q > 1$ instead $p = q = 1$ [Equation (2) in the next Section] within a "family of mettalic means" (Spinadel, 1998, 1999).

## 2. BASIC CONCEPTS

The GM arose from the division of a unit segment line AB into two parts (Fig. 1b): first $x$ and second $1 - x$, such that

$$\frac{x}{1-x}=\frac{1}{x}\ . \qquad (1)$$

On the other hand, one can say that GM follows from a square equation

$$x^2 \pm px - q = 0\,, \qquad (2)$$

where $p$ = 1, $q$ = 1, which solutions are:

$$x_{1,2}=\frac{-1\pm\sqrt{5}}{2}\,,\ \text{ or }\ x_{1,2}=\frac{1\pm\sqrt{5}}{2}\,. \qquad (3)$$

Stakhov (Stakhov, 1989, Figure 7 and Equation 26) revealed a possible generalization of GM, from which it follows:

$$x^n + x = 1, \qquad (4)$$

where $n$ = 1, 2, 3, ... .

## 3. A NEW GENERALIZATION

In this paper we reveal, however, a further generalization, such that in equation (2) $p$ = 1 and $q$ = $m$/2 ; thus, we consider the following equation:

$$x^n + x = \frac{m}{2}\ , \qquad (5)$$

where $n$ = 1, 2, 3, ... ; and $m$ = 0, 1, 2, 3, ... .

For $n$ = 2 and $m$ = 2, we have the well-known GM (Fig. 1b), and for other values – the generalized GM's (Fig. 1a, c, d and Table 1). By this, the values of $m$ correspond to the square roots of odd positive integers ($r$ = 1, 3, 5, 7, ...), through the generalized formula (3) as:

$$x_{1,2}=\frac{-1\pm\sqrt{r}}{2}\,,\ \text{ or }\ x_{1,2}=\frac{1\pm\sqrt{r}}{2}\,. \qquad (6)$$

In the following Fig. 1 we shall give the geometric and algebraic interpretations for $m$ = 1, 2, 3, 4.

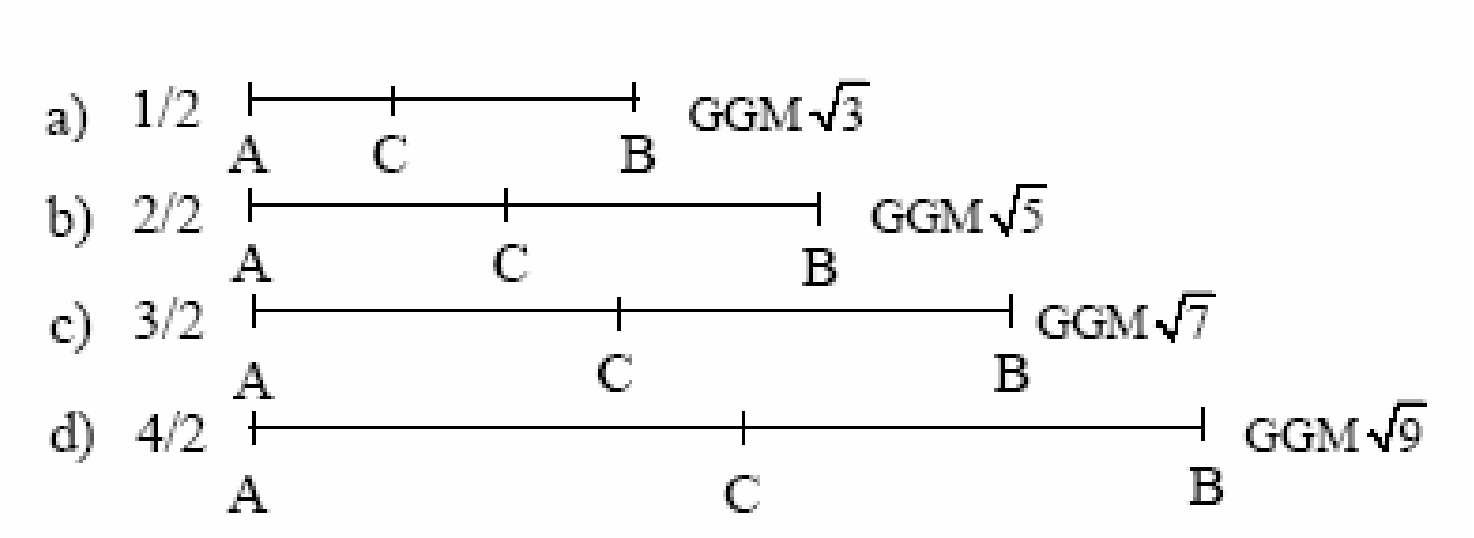

**Figure 1.** Generalized Golden mean by equation (5)

Consequently, it makes sense to speak about "Golden mean per root" of 1, of 3, of 5, of 7, of 9, and so on, respectively. Also, it makes sense to see the GM as an example of "the symmetry in the simplest case" (Marcus, 1989), just in the case when $n = 1$, $m = 2$. (Notice that this case is equivalent with the case $n = 2$, $m = 4$ as it is self-evident from figure 1d.)

### 3.1. Integer and non-integer solutions

In the following scheme (Table 1) we shall give the integer and non-integer solutions of Generalized GM.

From Table 1 it is evident that the sum of absolute values of solutions $x_1$ and $x_2$, to equation (5) equals, is $\sqrt{r}$, which represents the first cathetus (first leg) of triangle, $\sqrt{r} - m - h$. In such a triangle $m$ is the second cathetus and $h$ the hypotenuse. All such triangles on the left side in Table 1 appear as Diophantus' (Pythagorean) triangles (*see* Box 1), and on the right side their corresponding triangles. According to equations (2) and (6) there are four solutions, two positive and two negative, with two absolute values, as it is given in Table 1. Notice that $r - m - h$ triplets on the right side in Table 1 correspond to the Fibonacci triplets in first three cases (with ***h*** as an ordinal number)(Mišić, 2004): 0−1−1, 1−2−3, 2−3−5 through a growth for Fibonacci distance triplet 1−1−2. In next (forth) step, with the same distance 1−1−2, the Lucas' triplet 3−4−7 appears, which grows in all further steps just for one Fibonacci distance triplet 1−1−2. Notice also the next relations: on the left side in Table 1 the left-***h***, as well as the left-*m*, grows for 4k units (k = 0, 1, 2, 3, ...) whereas on the right side the right-***h*** and right-*m* grow just for one unit; the *r* on the left corresponds with $r^2$ on the right; the left-N$^{th}$ triangle appears in the right sequence through this "4k" regularity. (*Remark* 1: From the "4k" regularity follow triangles 0$^{th}$ – **0**$^{th}$, $1^{st} - \mathbf{4}^{th}$ ($1^{st}$ on the left, and $\mathbf{4}^{th}$ on the right in table 1), $2^{nd} - \mathbf{12}^{th}$, $3^{rd} - \mathbf{24}^{th}$, $4^{th} - \mathbf{40}^{th}$ etc., with next solutions: $[0+(4\times 0)=\mathbf{0}]$, $[0+(4\times 1)=\mathbf{4}]$, $[4+(4\times 2)=\mathbf{12}]$, $[12+(4\times 3)=\mathbf{24}]$, $[24+(4\times 4)=\mathbf{40}]$, etc.)

| N | $x_1$ | | $x_2$ | | $h$ | $m$ | $\sqrt{r}$ | | N | $x_1$ | | $x_2$ | | $h$ | $m$ | $\sqrt{r}$ |
|---|---|---|---|---|---|---|---|---|---|---|---|---|---|---|---|---|
| 0. | $0^2$ | + | $1^2$ | = | 1 | 0 | $\sqrt{1}$ | | 0. | $0^2$ | + | $1^2$ | = | 1 | 0 | $\sqrt{1}$ |
| | (0 | + | $1)^2$ | = | 1 | | | | | (0 | + | $1)^2$ | = | 1 | | |
| 1. | $1^2$ | + | $2^2$ | = | 5 | 4 | $\sqrt{9}$ | | 1. | $(x_1)^2$ | + | $(x_2)^2$ | = | 2 | 1 | $\sqrt{3}$ |
| | (1 | + | $2)^2$ | = | 9 | | | | | $(x_1$ | + | $x_2)^2$ | = | 3 | | |
| 2. | $2^2$ | + | $3^2$ | = | 13 | 12 | $\sqrt{25}$ | | 2. | $(x_1)^2$ | + | $(x_2)^2$ | = | 3 | 2 | $\sqrt{5}$ |
| | (2 | + | $3)^2$ | = | 25 | | | | | $(x_1$ | + | $x_2)^2$ | = | 5 | | |
| 3. | $3^2$ | + | $4^2$ | = | 25 | 24 | $\sqrt{49}$ | | 3. | $(x_1)^2$ | + | $(x_2)^2$ | = | 4 | 3 | $\sqrt{7}$ |
| | (3 | + | $4)^2$ | = | 49 | | | | | $(x_1$ | + | $x_2)^2$ | = | 7 | | |
| 4. | $4^2$ | + | $5^2$ | = | 41 | 40 | $\sqrt{81}$ | | 4. | $1^2$ | + | $2^2$ | = | 5 | 4 | $\sqrt{9}$ |
| | (4 | + | $5)^2$ | = | 81 | | | | | (1 | + | $2)^2$ | = | 9 | | |
| 5. | $5^2$ | + | $6^2$ | = | 61 | 60 | $\sqrt{121}$ | | 5. | $(x_1)^2$ | + | $(x_2)^2$ | = | 6 | 5 | $\sqrt{11}$ |
| | (5 | + | $6)^2$ | = | 121 | | | | | $(x_1$ | + | $x_2)^2$ | = | 11 | | |
| | (...) | | | | | | | | | (...) | | | | | | |

**Table 1.** The integer and non-integer solutions of Generalized Golden Mean

In the following example we shall consider the cases $n = 2$ and $m = 1, 2, 3$.

In the first case with $n = 2$ and $m = 1$, we have q = 0.5, the first case in Table 1 on the right (the first, not the zeroth), and in Figure 1a, as the case of "GM" per $\sqrt{3}$ with the two solutions given by equations (2) and (5):

$$x_1 = \left(-1+\sqrt{3}\right)/2 = 0.3660254...$$

and

$$x_2 = \left(-1-\sqrt{3}\right)/2 = -1.3660254... .$$

The satisfactory solution is positive solution $x_1$.

In the second case with $n = 2$, and $m = 2$ we have $q = 1$, the second case in Table 1 on the right, and in Figure 1b, as the case of GM per $\sqrt{5}$ with the two solutions:

$$x_1 = \left(-1+\sqrt{5}\right)/2 = 0.6180339...$$

and

$$x_2 = \left(-1-\sqrt{5}\right)/2 = -1.6180339... \; .$$

In the third case with $n = 2$ and $m = 3$, we have $q = 1.5$, as the case of "GM" per $\sqrt{7}$ , with the two solutions:

$$x_1 = \left(-1+\sqrt{7}\right)/2 = 0.8228756...$$

and

$$x_2 = \left(-1-\sqrt{7}\right)/2 = -1.8228756... \text{ , etc.}$$

## 4. THE METALLIC MEANS FAMILY

As it is known, it is very easy to find the members of "the metallic means family" (MMF) (Spinadel, 1999) as solutions of the equation (2), for various values of the parameters $p$ and $q$. In fact, if $p = q = 1$, we have the GM. Analogously, for $p = 2$ and $q = 1$ we obtain the Silver mean; for $p = 3$ and $q = 1$, we get the Bronze mean. For $p = 4$; $q = 1$ we have the next metallic mean, etc. On the other hand, if $p = 1$ and $q = 2$, we obtain the Copper mean. If $p = 1$ and $q = 3$, we get the Nickel mean and so on. Thus, we obtain all members of the MMF, which follow from square equation (2).

**Box 1. Diophantus' triangles**

$1^2 = 0^2 + 1^2$

$5^2 = 4^2 + 3^2$

$13^2 = 12^2 + 5^2$

$25^2 = 24^2 + 7^2$

$41^2 = 40^2 + 9^2$

$61^2 = 60^2 + 11^2$

$85^2 = 84^2 + 13^2$

*Diophantus' triangles are in fact the Pythagorean triangles with first cathetus "a" taken from the od natural number series (r = 1, 3, 5, 7, 9, ...), secon cathetus "b" calculated through the "4k regularity" presented in Remark 1, within Section 3.1, and hypotenuse "c", where c = b+1.*

**Box 2. The first Luca Pacioli's triangle**

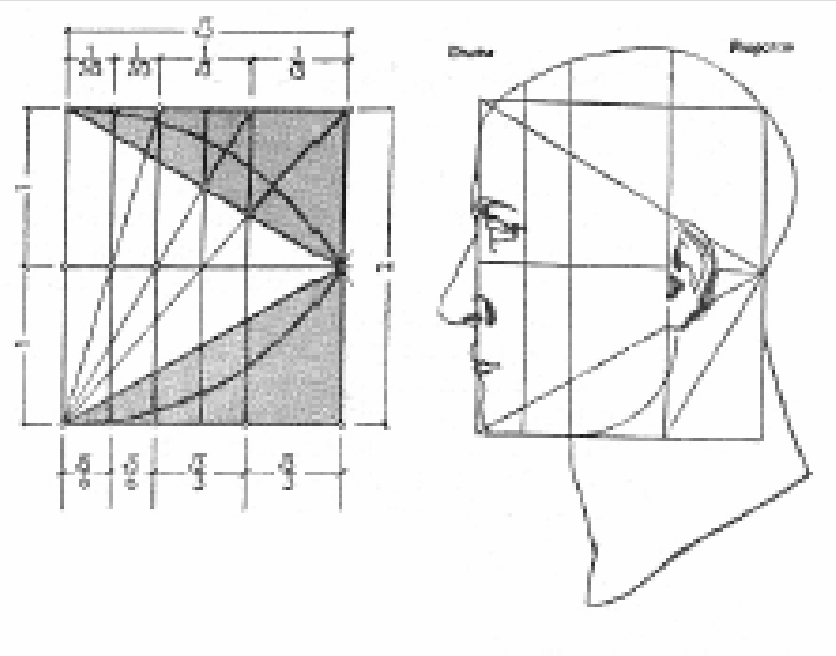

*This Luca Pacioli's triangle appears as the first triangle in the right side of Table 1. The right side of Figure follows from original Pacioli's picture, while on the left side there is a Zloković's mathematical analysis (Zloković, 1955).*

However, if we by (2) and (5) form the follow equation

$$x^n \pm px = \frac{m}{2}, \qquad (7)$$

where $n = 1, 2, 3, \ldots$ and $p = 1, 2, 3, \ldots$ , then we have a generalization of MMF; furthermore, we have a unification of "vertical" and "horizontal" generalization of GM.

Observe that De Spinadel (1999) found "the integer metallic means", for $q$ = 2, 6, 12, 20, 30, ..., which solutions ($x_1$, $x_2$), given by equation (2), are positive integers: (1, 2), (2, 3), (3, 4), (4, 5), (5, 6) ... (Spinadel, 1999, Section 3: "Furthermore, it is very easy to verify that ... the integer metallic means, $[2,\overline{0}]$, $[3,\overline{0}]$, $[4,\overline{0}]$, ..., appear in quite a regular way")(cf. Tables 2-3).

| | | | | | | | | | |
|---|---|---|---|---|---|---|---|---|---|
| 0 | **0** | 0 | 0 | 0 | 0 | 0 | 0 | 0 | 0 |
| **0** | 1 | **2** | 3 | 4 | 5 | 6 | 7 | 8 | 9 |
| 0 | **2** | 4 | **6** | 8 | 10 | 12 | 14 | 16 | 18 |
| 0 | 3 | **6** | 9 | **12** | 15 | 18 | 21 | 24 | 27 |
| 0 | 4 | 8 | **12** | 16 | **20** | 24 | 28 | 32 | 36 |
| 0 | 5 | 10 | 15 | **20** | 25 | **30** | 35 | 40 | 45 |
| 0 | 6 | 12 | 18 | 24 | **30** | 36 | **42** | 48 | 54 |
| 0 | 7 | 14 | 21 | 28 | 35 | **42** | 49 | **56** | 63 |
| 0 | 8 | 16 | 24 | 32 | 40 | 48 | **56** | 64 | **72** |
| 0 | 9 | 18 | 27 | 36 | 45 | 54 | 63 | **72** | 81 |

| | | |
|---|---|---|
| (0 x 0) + 0 = 00 | | 0 x 1 = 00 |
| (1 x 1) + 1 = 02 | | 1 x 2 = 02 |
| (2 x 2) + 2 = 06 | | 2 x 3 = 06 |
| (3 x 3) + 3 = 12 | | 3 x 4 = 12 |
| (4 x 4) + 4 = 20 | | 4 x 5 = 20 |
| (5 x 5) + 5 = 30 | | 5 x 6 = 30 |
| (6 x 6) + 6 = 42 | | 6 x 7 = 42 |
| (7 x 7) + 7 = 56 | | 7 x 8 = 56 |
| (8 x 8) + 8 = 72 | | 8 x 9 = 72 |
| (9 x 9) + 9 = 90 | | |

**Table 2 (left).** The harmonic multiplication Table of decimal numbering system. This Table contains "the integer metallic means", for $q$ = 0, 2, 6, 12, 20, 30, 42, 56 and 72 on the diagonal in form of doublets (pairs): 0-0, 2-2, 6-6, 12-12, 20-20, 30-30, 42-42, 56-56 and 72-72.
**Table 3 (right).** The key of the harmonic multiplication Table. This key is related to positive integers: (0, 1), (1, 2), (2, 3), (3, 4), (4, 5), (5, 6) which appear as solutions ($x_1$, $x_2$), given by equation (2); as solutions for above given $q$ ($q$ = 0, 2, 6, 12, 20, 30, ...).

From Table 1 it is self-evident that Spinadel's "integer metallic means" are related to the Diophantus' triangles too (*see* Box 1), as well as to the square roots of positive integers which are squares of odd integers; thus, $r$ = 1, 9, 25, 49, 81, 121, etc. On the other hand, the generalization, given by equation (5) is related to the square roots of all odd integers; thus, r = 1, 3, 5, 7, 9, 11, 13, etc.

## 5. THE EULER'S GENERALIZATION

In the history of mathematics it was known a conflict between the famous atheist philosopher Diderot and the famous religious mathematician Euler. ... One day Euler stepped for Diderot and stated: *"Sir,* $(a+b^n)/n = x$*, hence God exists; reply!"* (Eves, 1976). Diderot, as well as any one up to these days had no idea what Euler was talking about. We start here with the hypothesis (for further investigations) that Euler had the idea about "De divina proportione" of Luca Pacioli (1509) (*see* Box 2). However, after presented discussion in previous Sections of this paper we can suppose that this Euler's "divine equation" can be interpreted as a most possible generalization of GM for all cases discussed in this paper. Namely, in the case $a = b$, $n = 2$ and if and only if $x = 1/2$, we have, by equations (2) and (4) just the GM; moreover, GM stand then (accordingly to the principle "if one, then all") the case of one more extended generalization:

$x^n + x^{n-1} = 1$ (Stakhov, 1989, Equation 25), and/or of $x^n + px^{n-1} = m/2$, but that is the subject of a separate work.

## 6. THE GENERALIZATION THROUGH FIBONACCI SERIES

From the time Leonardo Fibonacci's book *Liber Abaci* (1202) was published, it is assumed that the golden mean can be generated starting from a series of the numbers which follow by addition two previous; in other words, by the Fibonacci series. On the other hand, from the time of Edouard Lucas (1842-1891), we have a generalization of this law by expanding it with the so-called Lucas' series. Now, however, we also present a generalization through sequences that (in correspondence with a series of natural numbers) follow after Lukas' series: the generalization with the post-Lucas series which are closely related to the Fibonacci series (Table 4).

| | | | | | | | | | |
|---|---|---|---|---|---|---|---|---|---|
| 0 | 1 | 1 | 2 | 3 | 5 | 8 | 13 | 21 | ... |
| | | 1 | 1 | 2 | 3 | 5 | 8 | 13 | |
| 1 | 1 | 2 | 3 | 5 | 8 | 13 | 21 | 34 | |
| | | 1 | 1 | 2 | 3 | 5 | 8 | 13 | |
| 2 | 1 | 3 | 4 | 7 | 11 | 18 | 29 | 47 | |
| | | 1 | 1 | 2 | 3 | 5 | 8 | 13 | |
| 3 | 1 | 4 | 5 | 9 | 14 | 23 | 37 | 60 | |
| | | 1 | 1 | 2 | 3 | 5 | 8 | 13 | |
| 4 | 1 | 5 | 6 | 11 | 17 | 28 | 45 | 73 | |
| | | 1 | 1 | 2 | 3 | 5 | 8 | 13 | |
| 5 | 1 | 6 | 7 | 13 | 20 | 33 | 53 | 86 | |
| ... | | | | | | | | | |

**Table 4.** The full generalization of Golden mean through Fibonacci sequences, in relation to natural numbers series.

## 7. CONCLUSION

As we can see from the discussion in previous five Sections, some known generalizations of Golden mean, and this new one, given here, appear to be the cases of one more extended generalization, given first by Luca Pacioly, and then by Leonhard Euler. On the other hand, bearing in mind that genetic code is determined by Golden mean (Rakočević, 1998b) one must takes answer to a question (in further researches), is that determination, or is not, valid for Generalized golden mean too. Certainly, the same question it takes set and for other natural and artificial systems.

## ACKNOWLEDGEMENTS

I am grateful to professors Ljubomir Ćirić, Djuro Koruga and Nataša Mišić for stimulating discussions and benevolent critique.

# Appendix

## Diagonal Fibonacci series

The diagonal Fibonacci series is related to the specific relation of Fibonacci numbers 2 and 3. Namely, 2/3 is the harmonic mean of a whole and its half, and 3/2 represents the limit of "golden numbers" (The Fibonacci Quarterly 1994, June-July, 211-217). Hence, it makes sense for these two numbers in the Fibonacci series to replace the place. This creates a new Fibonacci series, which in relation to Lucas series gives two new series in the form of diagonally coupled pairs of numbers (Table A1 and A2).

**Table A1.** Singlefold diagonal Fibonacci series

| | | | | | | | | | |
|---|---|---|---|---|---|---|---|---|---|
| 2 | 1 | 3 | 4 | 7 | 11 | 18 | 29 | 47 | 76 |
| | +2 | -1 | +1 | 1 x **0** | 1 x **1** | 1 x **1** | 1 x **2** | 1 x **3** | 1 x **5** |
| 3 | 2 | 5 | 7 | 12 | 19 | 31 | 50 | 81 | |
| Mirror symmetry: [**1** (001 / 100) **4**] | | | | | | | | | |
| | 76 | 123 | 199 | 322 | 521 | 843 | 1364 | 2207 | ... |
| | 1 x **5** | 1 x **8** | 1 x **13** | 1 x **21** | 1 x **34** | 1 x **55** | 1 x **89** | 1 x **144** | ... |
| 81 | 131 | 212 | 343 | 555 | 898 | 1453 | 2351 | 3804 | ... |

**Table A2.** Fourfold diagonal Fibonacci series

| | | | | | | | | | |
|---|---|---|---|---|---|---|---|---|---|
| 2 | 1 | 3 | 4 | 7 | 11 | 18 | 29 | 47 | 76 |
| | 4x0 | 4x1 | 4x1 | 4x2 | 4x3 | 4x5 | 4x8 | 4x13 | 4x21 |
| 3 | 2 | 5 | 7 | 12 | 19 | 31 | 50 | 81 | 131 |
| Mirror symmetry: [**4** (100 / 001) **1**] | | | | | | | | | |
| 47 | 76 | 123 | 199 | 322 | 521 | 843 | 1364 | 2207 | ... |
| | 4x21 | 4x34 | 4x55 | 4x89 | 4x144 | 4x233 | 4x377 | 4x610 | ... |
| | 131 | 212 | 343 | 555 | 898 | 1453 | 2351 | 3804 | ... |

From Figure A1 we see that tables A1 and A2 show the relations of Lucas' series and its first neighbor series. Precisely speaking in Tables A1 and A2 we have the Lucas' series versus its first neighbor series, with an arithmetic difference: single Fibonacci series versus fourfold Fibonacci series.

In the following tables are showed the relationships of further followers of the Lucas' series (Figure A1), such relationships as shown in the header of each table.

| | | | | | | | | | |
|---|---|---|---|---|---|---|---|---|---|
| (-1) | 1 | 0 | 1 | 1 | 2 | 3 | 5 | 8 | |
| | 1 | 1 | 2 | 3 | 5 | 8 | 13 | 21 | |
| (0) | **2** | **1** | **3** | **4** | **7** | **11** | **18** | **29** | |
| | 1 | 1 | 2 | 3 | 5 | 8 | 13 | 21 | |
| (1) | 3 | 2 | 5 | 7 | 12 | 19 | 31 | 50 | |
| | 1 | 1 | 2 | 3 | 5 | 8 | 13 | 21 | |
| (2) | 4 | 3 | 7 | 10 | 17 | 27 | 44 | 71 | |
| | 1 | 1 | 2 | 3 | 5 | 8 | 13 | 21 | |
| (3) | 5 | 4 | 9 | 13 | 22 | 35 | 57 | 92 | ... |
| ... | | | | | | | | | |

**Figure A1.** Lucas series (heavily shaded) with a previous neighbor and all followers. [Note: The "head" of Lukas' series consists of the first two numbers (2, 1) and the remaining numbers constitute the "body" of the series]

**Table A3.** The first series as the first adjacent to Lucas' series versus the second adjacent (follower), with an arithmetic difference: a double Fibonacci series versus a fivefold Fibonacci series (cf. Table A10)

| | | | | | | | | | |
|---|---|---|---|---|---|---|---|---|---|
| 3 | 2 | 5 | 7 | 12 | 19 | 31 | 50 | 81 | 131 |
| +2 | -2 | 2 x 0 | 2 x **1** | 2 x **1** | 2 x **2** | 2 x **3** | 2 x **5** | 2 x **8** | |
| 4 | 3 | 7 | 10 | 17 | 27 | 44 | 71 | 115 | ... |
| Mirror symmetry: [2 (010 / 101) 5] | | | | | | | | | |
| 3 | 2 | 5 | 7 | 12 | 19 | 31 | 50 | 81 | 131 |
| 5x**0** | 5x**1** | 5x**1** | 5x**2** | 5x**3** | 5x**5** | 5x**8** | 5x**13** | ... | |
| 4 | 3 | 7 | 10 | 17 | 27 | 44 | 71 | 115 | ... |

**Table A4.** Second series as the second adjacent to Lucas' series versus third adjacent, with arithmetic difference: irregular series versus sixfold Fibonacci series

| | | | | | | | | | |
|---|---|---|---|---|---|---|---|---|---|
| 4 | 3 | 7 | 10 | 17 | 27 | 44 | 71 | 115 | 186 |
| | -2 | 3 | 1 | 4 | 5 | 9 | 14 | 23 | 37 |
| 5 | 4 | 9 | 13 | 22 | 35 | 57 | 92 | 149 | ... |
| No mirror symmetry: [? (??? / 110) 6] | | | | | | | | | |
| 4 | 3 | 7 | 10 | 17 | 27 | 44 | 71 | 115 | 186 |
| | 6x**0** | 6x**1** | 6x**1** | 6x**2** | 6x**3** | 6x**5** | 6x**8** | 6x**13** | 6x**21** |
| 5 | 4 | 9 | 13 | 22 | 35 | 57 | 92 | 149 | 241 |
| | | | | | | | | | ... |

**Note**: The phrase "No mirror symmetry" means that the series is irregular, in the sense that its members are not the product of Fibonacci or Lukas numbers, nor are the product of the numbers from the series, adjacent to the sequence (followers) of the Lukas series.

**Table A5.** Third series as the third adjacent to Lucas' series versus fourth adjacent, with an arithmetic difference: double Lucas' series versus sevenfold Fibonacci series

| | | | | | | | | | |
|---|---|---|---|---|---|---|---|---|---|
| 5 | 4 | 9 | 13 | 22 | 35 | 57 | 92 | 149 | 241 |
| | -2 | 2 x **2** | 2 x **1** | 2 x **3** | 2 x **4** | 2 x **7** | 2 x **11** | 2 x **18** | 2 x **29** |
| 6 | 5 | 11 | 16 | 27 | 43 | 70 | 113 | 183 | ... |
| Mirror symmetry: [2 (0**1**0 / 1**1**1) 7] | | | | | | | | | |
| 5 | 4 | 9 | 13 | 22 | 35 | 57 | 92 | 149 | 241 |
| | 7x**0** | 7x**1** | 7x**1** | 7x**2** | 7x**3** | 7x**5** | 7x**8** | 7x**13** | 7x**21** |
| 6 | 5 | 11 | 16 | 27 | 43 | 70 | 113 | 183 | 296 |
| | | | | | | | | | ... |

**Table A6.** Fourth series as the fourth adjacent to Lucas' series versus the fifth adjacent, with an arithmetic difference: irregular series versus eightfold Fibonacci series

| | | | | | | | | | |
|---|---|---|---|---|---|---|---|---|---|
| 6 | 5 | 11 | 16 | 27 | 43 | 70 | 113 | 183 | 296 |
| | -2 | 5 | 3 | 8 | 11 | 19 | 30 | 49 | 79 ... |
| 7 | 6 | 13 | 19 | 32 | 51 | 83 | 134 | 217 | 351 |
| No mirror symmetry: [? (???? / 1000) 8] | | | | | | | | | |
| 6 | 5 | 11 | 16 | 27 | 43 | 70 | 113 | 183 | 296 |
| | 8x**0** | 8x**1** | 8x**1** | 8x**2** | 8x**3** | 8x**5** | 8x**8** | 8x**13** | 8x**21** |
| 7 | 6 | 13 | 19 | 32 | 51 | 83 | 134 | 217 | 351 |
| | | | | | | | | | ... |

**Table A7.** Fifth series as the fifth adjacent to Lucas' series versus sixth adjacent, with an arithmetic difference: double first adjacent series versus ninefold Fibonacci series

| | | | | | | | | | |
|---|---|---|---|---|---|---|---|---|---|
| 7 | 6 | 13 | 19 | 32 | 51 | 83 | 134 | 217 | 351 |
| | -2 | 2x**3** | 2x**2** | 2x**5** | 2x**7** | 2x**12** | 2x**19** | 2x**31** | 2x**50** ... |
| 8 | 7 | 15 | 22 | 37 | 59 | 96 | 155 | 251 | 406 |
| Mirror symmetry: [2 (001**0** /100**1**) 9] | | | | | | | | | |
| 7 | 6 | 13 | 19 | 32 | 51 | 83 | 134 | 217 | 351 |
| | 9x**0** | 9x**1** | 9x**1** | 9x**2** | 9x**3** | 9x**5** | 9x**8** | 9x**13** | 9x**21** |
| 8 | 7 | 15 | 22 | 37 | 59 | 96 | 155 | 251 | 406 |
| | | | | | | | | | ... |

**Table A8.** Sixth series as the sixth neighbor adjacent to the Luas' series versus the seventh neighbor, with an arithmetic difference: irregular series versus tenfold Fibonacci series

| | | | | | | | | | |
|---|---|---|---|---|---|---|---|---|---|
| 8 | 7 | 15 | 22 | 37 | 59 | 96 | 155 | 251 | 406 |
| | -2 | 7 | 5 | 12 | 17 | 29 | 46 | 75 | 121 ... |
| 9 | 8 | 17 | 25 | 42 | 67 | 109 | 176 | 285 | 461 |
| No mirror symmetry: [? (???? / 1010) 10] | | | | | | | | | |
| 8 | 7 | 15 | 22 | 37 | 59 | 96 | 155 | 251 | 406 |
| | 10x**0** | 10x**1** | 10x**1** | 10x**2** | 10x**3** | 10x**5** | 10x**8** | 10x**13** | 10x**21** |
| 9 | 8 | 17 | 25 | 42 | 67 | 109 | 176 | 285 | 461 |
| | | | | | | | | | ... |

**Table A9.** Seventh series as the seventh adjacent to Lukas' series versus the eighth adjacent, with an arithmetic difference: double second adjacent series versus eleventh Fibonacci sequence

| | | | | | | | | | |
|---|---|---|---|---|---|---|---|---|---|
| 9 | 8 | 17 | 25 | 42 | 67 | 109 | 176 | 285 | 461 |
| | -2 | 2x**4** | 2x**3** | 2x**7** | 2x**10** | 2x**17** | 2x**27** | 2x**44** | 2x**71** ... |
| 10 | 9 | 19 | 28 | 47 | 75 | 197 | 272 | 469 | 741 |
| Mirror symmetry: [2 (**0**010 / 101**1**) 11] | | | | | | | | | |
| 9 | 8 | 17 | 25 | 42 | 67 | 109 | 176 | 285 | 461 |
| | 11x**0** | 11x**1** | 11x**1** | 11x**2** | 11x**3** | 11x**5** | 11x**8** | 11x**13** | 11x**21** |
| 10 | 9 | 19 | 28 | 47 | 75 | 197 | 272 | 469 | 741 |

**Table A10.** Fibonacci series as a “neighbor” in itself: The original Fibonacci series versus the first adjacent to Lukas series

| | | | | | | | | | | | | | |
|---|---|---|---|---|---|---|---|---|---|---|---|---|---|
| **0** | **1** | **1** | 2 | 3 | 5 | 8 | 13 | 21 | 34 | 55 | 89 | 144 | 233 |
| | | | 1 | 1 | 0 | 1 | 1 | 2 | 3 | 5 | 8 | 13 | 21 |
| | | | 3 | 2 | 5 | 7 | 12 | 19 | 31 | 50 | 81 | 131 | 212 |
| | | | | | | | | | | | | | |
| 2 | 3 | 5 | 8 | 13 | 21 | 34 | 55 | 89 | 144 | 233 | 377 | 610 | 987 |
| | 3x**0** | 3x**1** | 3x**1** | 3x**2** | 3x**3** | 3x**5** | 3x**8** | 3x**13** | 3x**21** | 3x**34** | 3x**55** | 3x**89** | 3x**144** |
| 3 | 2 | 5 | 7 | 12 | 19 | 31 | 50 | 81 | 131 | 212 | 343 | 555 | 898 |
| Mirror symmetry: [2 (**0**10 / 01**1**) 3] | | | | | | | | | | | | | |
| 2 | 3 | 5 | 8 | 13 | 21 | 34 | 55 | 89 | 144 | 233 | 377 | 610 | 987 |
| | 2x**0** | 2x**1** | 2x**1** | 2x**2** | 2x**3** | 2x**5** | 2x**8** | 2x**13** | 2x**21** | 2x**34** | 2x**55** | 2x**89** | 2x**144** |
| 3 | 2 | 5 | 7 | 12 | 19 | 31 | 50 | 81 | 131 | 212 | 343 | 555 | 898 |

**Note 1**: The "head" of a Fibonacci series consists of the first three numbers (0, 1, 1), and the remaining numbers constitute the "body" of the series. Notice also the mirror symmetry of the head in the original and in the arithmetic difference series (green highlighted; cf the Note in legend of Figure A1). **Note 2**: The arrangement above the shaded area is given according to the idea of Prof. Djuro Koruga (private communication), and I am very grateful to him for that.